# Galois theory without abstract algebra

Leonid Lerner*

School of Physical Sciences, Flinders University, Adelaide, Australia 5001

**Abstract**

Galois theory is developed using elementary polynomial and group algebra. The method follows closely the original prescription of Galois, and has the benefit of making the theory accessible to a wide audience. The theory is illustrated by a solution in radicals of lower degree polynomials, and the standard result of the insolubility in radicals of the general quintic and above. This is augmented by the presentation of a general solution in radicals for all polynomials when such exist, and illustrated with specific cases. A method for computing the Galois group and establishing whether a radical solution exists is also presented.

## *1. Introduction*

One is frequently required to solve polynomial equations in the natural sciences [1]. For example, the normal mode frequencies $\omega$ of coupled oscillators, such as the masses on a spring shown in Fig, 1a, are the roots of the matrix eigenvalue equation

*Email: leonid.lerner@flinders.edu.au

$$\begin{vmatrix} m\omega^2 - 2k & k & 0 & .. \\ k & m\omega^2 - 2k & k & .. \\ 0 & k & m\omega^2 - 2k & .. \\ .. & .. & .. & .. \end{vmatrix} = 0. \qquad (1)$$

Similarly, coupled mode equations of optical waveguides, or a multiple level quantum system, Fig. 1b, reduce to polynomials whose roots are propagation constants and energy eigenvalues respectively

$$\begin{vmatrix} E - E_1 & C_{12} & C_{13} & .. \\ C_{12}^* & E - E_2 & C_{23} & .. \\ C_{13}^* & C_{23}^* & E - E_3 & .. \\ .. & .. & .. & .. \end{vmatrix} = 0. \qquad (2)$$

For constant coupling, Eqs. (1-2) are polynomial equations with degree determined by the number of coupled systems. It is well known that for completely arbitrary coefficients solutions in radicals, meaning expressions of the type shown in Table 1, exist only for polynomials of degree less than 5.

| # | Expression | | | |
|---|---|---|---|---|
| 1 | $A^a + B$ | | | |
| 2 | $A^a + B^b + C$ | $(A^a + B)^b + C$ | | |
| 3 | $A^a + B^b + C^c + D$ | $(A^a + B^b + C)^c + D$ | $((A^a + B)^b + C)^c + D$ | $(A^a + B)^b + C^c + D$ |

**Table 1**. Types of radical expressions in rational functions of polynomial coefficients $A, B, C, D$ listed in increasing number of root extractions.

For higher degrees, radical solutions exist only in specific cases, but there are some general solutions in special functions [**2**]. The reason for this lies in 'Galois Theory', but the details are mostly unfamiliar to those outside of pure mathematics due to the theory's extensive use of abstract algebra. Indeed, it is generally acknowledged that the intricate structure of the modern theory makes it hard to grasp in its entirety, and it is therefore easily forgotten [**3**].

It may therefore come as a surprise that the work of Galois, published in a posthumous memoir in 1846 [**4**], is almost entirely devoid of the abstractness of modern Galois Theory. Indeed, the bulk of the modern approach is based largely on later work [**5**], which starting from the original framework of Galois, put the theory on a more rigorous footing, and widened its scope. For this reason current elementary presentations of Galois theory at the technical level [**5-8**], use an extensive amount of abstract algebra. The latter, however, is unfamiliar to a wide audience outside of pure mathematics, which is more accustomed to algebraic operations on polynomials and group elements, and thus potentially able to appreciate the theory based on the original approach of Galois. Galois' memoir [**4, 8**], despite its ultimate simplicity, is quite hard to follow, due to some gaps in the presentation (for which it was rejected by the Paris Academy in 1830). However by filling in the gaps a fairly complete picture emerges (some of this has already been demonstrated in [**8**]). The aim of this paper is therefore to provide an exposition of Galois Theory, based on the original approach of Galois, which can act as an introduction for those who want to appreciate the theory at a technical level with the minimum effort. It goes without saying, that the simplicity and conciseness of the present treatment comes at the cost of some rigour in places, compared to modern Galois Theory. However such omissions can be made rigorous. The interested reader is referred to the many excellent introductory textbooks on the subject [5-8].

The essence of Galois' theory lies in a procedure it gives for associating a given polynomial equation with a group corresponding to the symmetry of its roots (the Galois group). It turns out that solutions expressible in radicals, examples of which are shown in Table 1, satisfy polynomials with groups of a particular type. One can then investigate the group of a given equation to see if it belongs to this type. If it does, Galois theory provides a method for solving it, if it does not, no solution in radicals can be obtained at all. The Galois group of an *n*-degree polynomial with arbitrary coefficients is the permutation group of $n$ objects, denoted $S_n$. Groups $S_5$ and above include permutations which do not belong to groups of equations solvable in radicals, explaining the insolubility of arbitrary polynomials of degree five and above.

The outline of this paper is as follows. At the outset elementary symmetric polynomials and the theorem of symmetric polynomials is introduced, providing a natural entry to the Galois Group. A unified way to solve $2^{nd}$, $3^{rd}$, and $4^{th}$ degree polynomials is then developed bearing a close correspondence to these group ideas. Since the method breaks down at $5^{th}$ degree polynomials the question arises as to whether another method exists. To this end the concept of an irreducible polynomial is introduced and some of its properties established. It is then demonstrated that a solution in radicals entails the splitting of the Galois group of an irreducible polynomial into a series of normal subgroups by the 'adjunction' of solutions of $x^p = A$. Moreover such a splitting corresponds exactly to the method by which the lower degree polynomials were solved, and therefore if a radical solution exists, it must proceed along the same lines. It is then shown that $S_5$, the Galois group of the general 5th degree polynomial, can not be split in such a way and so the general 5th degree polynomial in insoluble. An example of a solvable higher degree polynomial is then presented, its Galois group is calculated and shown to split in the required fashion. Galois theory is then used to derive its radical solution.

Galois' work was preceded by the Abel-Ruffini proof of the insolubility in radicals of 5th degree polynomials [**5, 8**]. This proof is both more elaborate and more restrictive than Galois Theory, giving no indication when a particular polynomial is solvable in radicals, or providing a means to solve it, hence it will not be discussed here.

## *2. Radical solutions of lower degree polynomials*

### **2.1 Elementary symmetric polynomials**

**Definition 2.1.** *The elementary symmetric polynomials $\sigma_i(r_1..r_n)$ in n variables are obtained by equating coefficients in the linear factorization and the power expansion of a monic polynomial of degree n .*

Thus for 2nd and 3rd degree polynomials we have

$$(x+r_1)(x+r_2) = x^2 + \sigma_1 x + \sigma_2$$
$$\sigma_1 \equiv r_1 + r_2$$
$$\sigma_2 \equiv r_1 r_2 \qquad (3)$$

$$(x+r_1)(x+r_2)(x+r_3) = x^3 + \sigma_1 x^2 + \sigma_2 x + \sigma_3$$
$$\sigma_1 \equiv r_1 + r_2 + r_3$$
$$\sigma_2 \equiv r_1 r_2 + r_1 r_3 + r_2 r_3 \qquad (4)$$
$$\sigma_3 \equiv r_1 r_2 r_3$$

Extending this to the general case we have

$$\begin{aligned}\sigma_1 &\equiv \sum_i r_i \\ \sigma_2 &\equiv \sum_{i<j} r_i r_j \\ \sigma_3 &\equiv \sum_{i<j<k} r_i r_j r_k \\ &\ldots\end{aligned} \qquad (5)$$

It follows that all $\sigma_i(r_1..r_n)$ can be rationally expressed in terms of $\sigma_i(r_1..r_{n-1})$ and the extra variable $r_n$. Denoting $\sigma_i(r_1..r_{n-1})$ by $\tau_1..\tau_{n-1}$ for clarity, Eq. (5) gives

$$\begin{aligned}\sigma_1 &= r_n + \tau_1 \\ \sigma_2 &= r_n \tau_1 + \tau_2 \\ \sigma_3 &= r_n \tau_2 + \tau_3 \\ &\quad .. \\ \sigma_n &= r_n \tau_{n-1}\end{aligned} \qquad (6)$$

The reverse is also true, and $\tau_1..\tau_{n-1}$ can be rationally expressed in terms of $\sigma_1..\sigma_n$ and $r_n$. Starting from the first line of Eq. (6), $\tau_1 = \sigma_1 - r_n$, and this can be substituted into the second relation to give $\tau_2$, and so on.

## 2.2 Fundamental theorem of symmetric polynomials

**Proposition 2.1** *Every symmetric polynomial $f(r_1..r_n)$ can be rationally expressed as a polynomial in the elementary symmetric polynomials, $f(r_1..r_n) = g(\sigma_1..\sigma_n)$.*

There is a simple proof by induction. Assume the proposition is valid for polynomials in $n-1$ variables, i.e. there is a polynomial $g$ such that $f(r_1..r_{n-1}) = g(\sigma_1..\sigma_{n-1})$. Setting $r_n = 0$ turns a symmetric polynomial in $n$ variables into one in $n-1$ variables. The assumption is therefore equivalent to

$$[f(r_1..r_n) - g(\sigma_1..\sigma_n)]_{r_n=0} = 0. \tag{7}$$

But this means that $r_n$ is a root of $f - g$ and hence divides $f - g$. Then by symmetry so do all the other $r_i$

$$\frac{f(r_1..r_n) - g(\sigma_1..\sigma_n)}{r_1..r_n} \equiv f^1(r_1..r_n). \tag{8}$$

The new polynomial $f^1$ is clearly symmetric, and of lower degree than $f$. We can now repeat the procedure of Eqs. (7-8) on $f^1$, lowering the degree until at some iteration we must achieve $f^n = g$ since there are no symmetric polynomials of lower degree then the elementary symmetric polynomials. The theorem is thus proved.

*Example*

$$f(r_1, r_2) \equiv r_1^4 + r_2^4, \quad \sigma_1 \equiv r_1 + r_2, \quad \sigma_2 \equiv r_1 r_2.$$

From Eq. (7)

$$\left[r_1^4 + r_2^4 - (r_1 + r_2)^4\right]_{r_2=0} = 0. \tag{9}$$

Applying Eq. (8)

$$\frac{r_1^4 + r_2^4 - (r_1 + r_2)^4}{r_1 r_2} = -4(r_1^2 + r_2^2) - 10 r_1 r_2. \tag{10}$$

Repeating this once more we achieve the required decomposition

$$r_1^4 + r_2^4 = \sigma_1^4 - 4\sigma_1^2 \sigma_2 + 2\sigma_2^2. \tag{11}$$

## 2.3 The Galois Group

In the previous section we showed that any function of the roots of a polynomial invariant under all permutations of these roots is 'known' in terms of the polynomial coefficients, in the sense that Eqs. (7-8) provide a procedure for calculating it. We can generalise this idea in the following definition

**Definition 2.2.** *The Galois group of a polynomial is the maximum group of permutations under which any function of its roots must be invariant, for it to be a known function of the polynomial coefficients.*

The group property arises because if the result of two permutations, say $r_1 \leftrightarrow r_2$ and $r_2 \leftrightarrow r_3$ is known, then so is the result of a combined permutation $r_1 \leftrightarrow r_3$. The previous section showed that if the polynomial coefficients are completely arbitrary, invariance under all permutations of the roots is sufficient for a function to be known. This is also necessary because asymmetric polynomials can not be constructed from rational combinations of symmetric functions unless some other relations exist between the roots.

The Galois group of an arbitrary polynomial of degree $n$, is therefore the permutation group of $n$ objects, $S_n$. Making the polynomial less general, for instance by specifying the coefficients such that there exist *rational* relations between its roots, may make the Galois group smaller. A method for calculating the Galois group is given in Section 4.

## 2.4 Quadratic polynomials

Best known are expressions for the roots of quadratic equations, obtained by 'completing the square'. Here it is more convenient to think in terms of symmetric polynomials. Thus the functions

$$\sigma_1 \equiv r_1 + r_2, \quad t^2 \equiv (r_1 - r_2)^2, \tag{12}$$

are symmetric under all permutations of the roots (the Galois group of the quadratic is $S_2$), and are therefore known. Moreover, the roots $-r_1, -r_2$ are linearly related to these symmetric functions $\sigma_1$ and $t$ and can therefore be easily solved

$$r_1, r_2 = (\sigma_1 \pm t)/2. \tag{13}$$

Rewriting $t^2$ in terms of the elementary symmetric polynomials we have

$$r_1, r_2 = \left(\sigma_1 \pm (\sigma_1^2 - 4\sigma_2)^{1/2}\right)/2, \tag{14}$$

and we have solved the quadratic equation in radicals.

## 2.5 Cubic polynomials

The key to solving the quadratic was finding a power of a linear combination of the roots invariant under $S_2$, called a 'resolvent'. In Eq. (12) the coefficients in the resolvent are the square roots of unity, $\alpha^2 = 1$. For the cubic $x^3 + ax^2 + bx + c = 0$, it is therefore natural to consider the cube roots of unity $\alpha_1 = 1$, $\alpha_2 = e^{2\pi i/3}$, $\alpha_3 = e^{4\pi i/3}$.

By reference to Eq. (12) we have

$$\sigma_1 = r_1 + r_2 + r_3, \quad t_1^3 \equiv (r_1 + \alpha_2 r_2 + \alpha_3 r_3)^3, \quad t_2^3 \equiv (r_1 + \alpha_2 r_3 + \alpha_3 r_2)^3 \tag{15}$$

Since cyclically permuting the roots $r_i$ in either of the expressions for $t_1^3$ or $t_2^3$ multiplies $t_1$ or $t_2$ by a cube root of unity, both $t_1^3$ and $t_2^3$ are invariant under the permutations of the cyclic group, $Z_3$. However unlike the quadratic, this is insufficient for them to be known, because $S_3$ has two cyclic subgroups:

$$
\begin{array}{ccc}
& S_3 & \\
Z_3 & & Z_3' \\
\begin{pmatrix} 1,2,3 \\ 2,3,1 \\ 3,1,2 \end{pmatrix} & \leftarrow g = (1,2) \rightarrow & \begin{pmatrix} 2,1,3 \\ 1,3,2 \\ 3,2,1 \end{pmatrix}
\end{array}
\tag{16}
$$

and a known quantity must be invariant under both sets of cycles.

The symmetry of $t_1^3$ and $t_2^3$ is intermediate between that of the roots (no symmetry), and $S_3$ required for a function to be known. We can put this increased symmetry to use by observing that symmetric combinations of $t_1^3$ and $t_2^3$ are invariant under permutations of both $Z_3$ and $Z_3'$, and are therefore known. In particular, we can calculate $t_1^3 + t_2^3$ and $(t_1^3 - t_2^3)^2$ in terms of the roots using Eq. (15), and then use Eqs. (7-8) to relate these to the cubic coefficients

$$t_1^3 + t_2^3 = 9ab - 2a^3 - 27c, \quad (t_1^3 - t_2^3)^2 = 81(a^2 b^2 - 4a^3 c + 18abc - 27c^2 - 4b^3). \tag{17}$$

Equation (17) is easily solved for $t_1^3$ and $t_2^3$. Extracting the cube roots gives us the resolvents $t_1$ and $t_2$, which are linearly related to the roots through Eq. (15). The cubic is thus solved.

Group $S_3$ consists of two presentations (or cosets) of the cyclic subgroup $Z_3$ and $Z_3'$, linked by the permutation $g \equiv (1,2)$. The permutation $n \equiv (1,2,3)$ applied to any element of either coset generates all its other elements, so that $g$ is a 'generator' of $Z_3$, which is thus Abelian. The elements of $Z_3'$ in Eq. (16) are ordered by $Z_3 \xrightarrow{g} Z_3'$, which differs from the order in which they are generated by $n \equiv (1,2,3)$, so $g$ and $n$ do not commute, indeed $gn^2 = ng$, and thus $S_3$ is 'non-Abelian'.

A subgroup $N$ of group $G$ is normal if conjugating any $n_i \in N$ with all $g \in G$ returns an element of $N$, i.e. $g n_i g^{-1} = n_j$ for some $j$. Put another way, the left and right cosets $gN$ and $Ng$ to which $gn_i$ and $n_i g$ belong, must be identical. Coset $Ng$ obtained by applying the elements of subgroup $N$ to a starting element $g$ is a presentation of the subgroup, that is the permutations linking the elements of this coset are those of $N$, however coset $gN$, obtained by applying a single permutation $g$ to all the elements of $N$ is not so, unless $N$ is normal. In Eq. (16), $gZ_3 = Z_3'$ is clearly a presentation of $Z_3$, thus $Z_3$ is a normal subgroup of $S_3$. The permutations linking the cosets of a normal subgroup form a 'quotient' group. That is, the product of two elements $a_i, b_i$ in cosets $An$ and $Bn$, with $A, B \in Q$ belong to the same coset $ABn$, viz $a_i b_j = An_i Bn_j = ABn_i' n_j$, so $Q$ forms a group in this sense. The quotient group $Q$ is denoted $G/N = Q$. In this case it is $Z_2$, so that $S_3 / Z_3 = Z_2$, or $1 \triangleleft Z_3 \stackrel{Z_2}{\triangleleft} S_3$.

## 2.6 Quartic polynomials

Since the composition of $S_3$ into normal subgroups solved the cubic, to solve the quartic $x^4 + a_3 x^3 + a_2 x^2 + a_1 x + a_0 = 0$ we consider the composition series [6] of $S_4$

$$1 \stackrel{Z_2}{\triangleleft} Z_2 \stackrel{Z_2}{\triangleleft} D_2 \stackrel{Z_3}{\triangleleft} A_4 \stackrel{Z_2}{\triangleleft} S_4, \tag{18}$$

see Fig. 2. The first normal subgroup in the series from $S_4$ is the 'alternating subgroup' $A_4$, which is the set of all even permutations (consisting of an even number of transpositions). It forms a group because a sequence of even transpositions is even. It is normal, because it is related to its other presentation $A_4'$ by a single permutation, say (1,2). The solution now follows the same procedure as for the cubic. We find a function symmetric under $A_4$, say $\theta_1$. Let $\theta_2$ be its image under the permutation (1,2) (which belongs to the $S_4 / A_4 = Z_2$ quotient group). Then the following combinations exhibit $S_4$ symmetry and are thus known in terms of the polynomial coefficients $a_i$

$$\theta_1 + \theta_2 = f(a_i) \quad (\theta_1 - \theta_2)^2 = g(a_i). \tag{19}$$

Now we find a function symmetric under $D_2$, say $t_1$. Let $t_2$ be its image under the permutation (2,3,4) (which belongs to the $A_4 / D_2 = Z_3$ quotient group). Since $Z_3$ is a 3-cycle there is a third image $t_3$. Then the following combinations exhibit $A_4$ symmetry and are thus known in terms of $\theta_i$

$$(t_1 + \alpha_1 t_2 + \alpha_2 t_3)^3 = \theta_1, \quad (t_2 + \alpha_1 t_1 + \alpha_2 t_3)^3 = \theta_2, \quad (t_1 + \alpha_1 t_3 + \alpha_2 t_2)^3 = \theta_3. \tag{20}$$

To complete the solution we split $D_2$ into its normal subgroups

$$
\begin{array}{ccc}
& D_2 & \\
Z_2 & & Z_2' \\
\begin{pmatrix} 1,2,3,4 \\ 2,1,4,3 \end{pmatrix} & \leftarrow (1,3)(2,4) \rightarrow & \begin{pmatrix} 3,4,1,2 \\ 4,3,2,1 \end{pmatrix} \\
d_1 = r_1 + r_2 - r_3 - r_4 & & d_2 = r_3 + r_4 - r_1 - r_2
\end{array}
\tag{21}
$$

Then the following combination exhibits $D_2$ symmetry and is thus known in terms of $t_i$

$$(d_1^2 + a_3^2)/4 = (r_1 + r_2)^2 + (r_3 + r_4)^2 = t_1 \tag{22}$$

$d_1 = r_1 + r_2 - r_3 - r_4$ is found by symmetrizing $z_1 = r_1 + ir_2 - r_3 - ir_4$, a linear combination of the roots invariant only under the identity, with respect to the permutations of $Z_2$ (these permutations satisfy the $Z_2$ multiplication table).

Combining Eq. (22) with $a_3 = r_1 + r_2 + r_3 + r_4$ gives a quadratic equation for $r_1 + r_2$. Repeating this for the other images of $D_2$ all $r_i + r_j$ can be found and the 4 roots determined. The quartic is thus solved.

### 2.7 General procedure for solvable polynomials

The previous section demonstrated on the example of the quartic, the procedure for finding the roots of a polynomial solvable in radicals. Starting with the polynomial coefficients exhibiting the symmetry of Galois group of the equation, we work down, lowering the symmetry by reducing the size of the Galois group, until at the bottom level, we relate these coefficient to the roots. The key to this procedure is the fact that a single permutation relates presentations of the intermediate groups, i.e. that the subgroups are normal with cyclic quotient groups. Therefore a polynomial with Galois group $Q$ and the composition series $1 \triangleleft ... C \overset{Z_c}{\triangleleft} B \overset{Z_b}{\triangleleft} A \overset{Z_a}{\triangleleft} Q$, has the solution

$$\begin{aligned} Q_i &= \left(A_{p_{1i}} + \alpha_a A_{p_{2i}} + ...\alpha_a^{a-1} A_{p_{ai}}\right)^a \\ A_i &= \left(B_{p_{1i}} + \alpha_b B_{p_{2i}} + ...\alpha_b^{b-1} B_{p_{bi}}\right)^b \\ B_i &= \left(C_{p_{1i}} + \alpha_c C_{p_{2i}} + ...\alpha_c^{c-1} C_{p_{ci}}\right)^c \\ &... \qquad\qquad ... \\ Z_i &= \phantom{(}r_{p_{1i}} + q_1 r_{p_{2i}} + ...q_{n-1} r_{p_{ni}} \end{aligned} \qquad (23)$$

where the top level functions $Q_i$ are expressible in terms of the equation coefficients, and the bottom level functions $Z_i$, are linear combination of the $n$ roots $r_i$, with the constants $q_i$ chosen so that $Z_i$ are invariant only under the identity. The functions $A_i, B_i..$ are invariant under permutations of groups $A, B..$ respectively, and are obtained from $A_1, B_1..$ by successive application of the single permutation (belonging to the quotient groups $Z_A, Z_B..$respectively) by which one passes from one presentation of the group to the next. The quantities $p_{ij}$ designates elements of the cyclic permutation, while $\alpha_i$ are the $i^{\text{th}}$ roots of unity.

## 2.8 Quintic polynomials

The solution of Section 2.7 does not apply to the general quintic because $S_5$ can not be composed in a series of normal subgroups with cyclic quotient groups. Therefore the quintic is either insoluble in radicals, or another solution method exists. Section 3.3 demonstrates that the only solution in radicals is that of Section 2.7.

To show that $S_5$ can not be composed as required consider the cyclic subgroup $Z_3$ generated by $(1,2,3)$, contained in $S_n$. Groups $S_5$ and higher contain two additional distinct symbols $4,5$, and one can write $(1,2,3) = g^{-1}h^{-1}gh$, where $g \equiv (4,2,1), h \equiv (1,5,3)$. If $N$ is a normal subgroup of $S_5$ the cosets to which $g,h$ belong are linked to $N$ by a single permutation of the quotient group, say $G,H$ respectively: $\{g\} \xleftarrow{G} \{n\} \xrightarrow{H} \{h\}$. Then

$$g^{-1}h^{-1}gh = n_g^{-1}G^{-1}n_h^{-1}H^{-1}Gn_g Hn_h \tag{24}$$

where $n_i$ are elements of $N$. Since the quotient group is cyclic, $G,H$ commute, and they can be shuffled past elements of $N$ until they combine to form the identity. The cycle $(1,2,3)$ is then a product of elements of $N$, and therefore also belongs to $N$. The entire series of 3-cycles is then passed inside $N$ from one step in the composition sequence to the next, where the argument can be repeated. The subgroup $N$ thus always contains all 3-cycles and can never be reduced to the identity. Then $S_n$ with $n \geq 5$ does not have the required composition series.

A solution in radicals still exists for restricted quintics whose Galois subgroup of $S_5$, has a normal composition. The largest such subgroup is the Frobenius group $F_{20}$, consisting of the cyclic permutations of (1,2,3,4,5) (1,3,5,2,4) (1,5,4,3,2) (1,4,2,5,3), that is 20 elements in all. The structure of this group is similar to $S_3$ of the cubic, with two generators $g \equiv (1,2,3,4,5) \in Z_5$ and $h \equiv (2,3,5,4)$ satisfying the same relation $hg^2 = gh$. One passes between the 4 presentations of the cyclic subgroup $Z_5$ by the permutation

(2,3,5,4), hence $Z_5$ is a normal subgroup. The composition series of $F_{20}$ is therefore

$1 \triangleleft Z_5 \overset{Z_2}{\triangleleft} D_5 \overset{Z_2}{\triangleleft} F_{20}$, and the Galois groups of solvable quintics must belong to this set.

A polynomial invariant under just the permutations of $F_{20}$ can be formed by symmetrizing $r_1^2 r_2 r_5$ under its 20 permutations, producing a 10 term sum

$$t_1 = r_1^2 r_2 r_5 + r_2^2 r_3 r_1 + r_3^2 r_4 r_2 + r_4^2 r_5 r_3 + r_5^2 r_1 r_4 + \leftrightarrow (2,3,5,4) \; . \tag{25}$$

Since $F_{20}$ contains 20 of the 120 permutations in $S_5$, there are six presentations of $F_{20}$ in $S_5$, corresponding to six invariants $t_i$. As $F_{20}$ consists of all 5 and 4 cycles these can be obtained by applying $F_{20}$ to the 3-cycle (3,4,5) and two cycle (4,5). Since $F_{20}$ is solvable, the five roots of the quintic $r_i$ can be solved [**10**] in terms of $t_i$ as described in the previous section. However since $S_5$ is insoluble, there is no radical solution for the $t_i$ in terms of the quintic coefficients (the $S_5$ invariants), and the best we can do is as follows.

Symmetric combinations of the six $t_i$, are $S_5$ invariant, and hence form a known sextic, $(t - t_1)..(t - t_6)$. For a quintic in Bring form: $x^5 + ax + b = 0$, to which all quintics can be reduced in radicals, this sextic is [**11**]

$$(t + 2a)^4 (t^2 + 16a^2) = 5^5 b^4 (t + 3a). \tag{26}$$

Since there is no group in $S_5$ with $F_{20}$ as normal subgroup this sextic can not be solved in radicals unless $t$ is rational. This condition can be parameterized as [**11**]

$$x^5 + 5d^4 \frac{4c - 3}{c^2 + 1} x + 4d^5 \frac{2c + 11}{c^2 + 1} = 0, \tag{27}$$

where $d$ is any real number but $c$ must be rational. For solvable quintics the Bring form is either that of Eq. (27), or can be factorised rationally into lower degree polynomials. Section 4 shows how to establish this by direct calculation.

## 3 When is a polynomial solvable in radicals?

The central result of Galois' theory states that a solution in radicals must necessarily be of the form of Section 2.7. To demonstrate this result we shall first recall some properties of irreducible polynomials.

### 3.1 Irreducible polynomials

**Definition 3.1.** *An irreducible polynomial is one that can not be factored into lower degree polynomials whose coefficients are rational functions of the original polynomial coefficients.*

Thus $f(x) \equiv x^2 - 2$ can not be rationally factored as $f(x) = g(x)h(x)$ with $f(x)$ and $g(x)$ being rational lower degree polynomials, and is thus irreducible.

**Definition 3.2.** *The 'adjunction' of a radical $\sqrt{A}$ is an extension of polynomial $g(x)$ to $g(x,\sqrt{A})$ such that the coefficients of $g(x)$ are rational functions of the new parameter.*

Thus $f(x) \equiv x^2 - 2$ can be factored into rational polynomials in $x$ and $\sqrt{2}$, since $x^2 - 2 = (x + \sqrt{2})(x - \sqrt{2})$. In this way a properly chosen adjunction can reduce an irreducible polynomial. Adjunction can also operate in stages as the following example illustrates.

Example:

$f(x) = x^4 - 10x^2 + 1$ is an irreducible quartic. Adjunction of $\sqrt{2}$ reduces it to a product of quadratic polynomials

$$g(x, \sqrt{2}) = x^2 - 2\sqrt{2}x - 1, \quad h(x, \sqrt{2}) = x^2 + 2\sqrt{2}x + 1. \tag{28}$$

Further adjunction of $\sqrt{3}$ completely factorizes it into a product of four linear terms

$$g(x, \sqrt{2}, \sqrt{3}) = x - \sqrt{2} - \sqrt{3}, \quad g'(x, \sqrt{2}, \sqrt{3}) = x - \sqrt{2} + \sqrt{3} \tag{29}$$

$$h(x, \sqrt{2}, \sqrt{3}) = x + \sqrt{2} - \sqrt{3}, \quad h'(x, \sqrt{2}, \sqrt{3}) = x + \sqrt{2} + \sqrt{3}$$

**Proposition 3.1.** *An irreducible polynomial with rational coefficients can not have any root in common with another distinct polynomial with rational coefficients without having all its roots in common. As a corollary, given an irreducible polynomial there exists no polynomial of lower degree that has any roots in common with that irreducible polynomial.*

This can be proved using Euclid's algorithm. Let $f(x)$ be an irreducible polynomial, and let $g(x)$ be another polynomial of lower degree, with which it shares roots $r_i$. Using polynomial division we can *rationally* determine a quotient polynomial $q(x)$ and remainder polynomial $g_1(x)$ which satisfy

$$f(x) = q(x)g(x) + g_1(x), \tag{30}$$

where the degree of the remainder $g_1(x)$ is less than the divisor $g(x)$. At the common roots $r_i$, both $f(x)$ and $g(x)$ are zero, consequently we must also have $g_1(r_i) = 0$. Hence $r_i$ are also the roots of the remainder polynomial. We can now repeat the process, dividing $g(x)$ by $g_1(x)$ to obtain a new remainder polynomial $g_2(x)$ of lower degree than $g_1(x)$ which also shares the common roots. Eventually to avoid a contradiction we must reach $g_n(x) = 0$ so that all the roots of $g_{n-1}(x)$ are common to $f(x)$ and $g(x)$. This means that $g_{n-1}(x)$ is a rationally determined factor polynomial of $f(x)$, and $f(x)$ is reducible, which contradicts the assumption.

**Proposition 3.2.** *An irreducible polynomial can not have repeated roots.*

Thus a polynomial with repeated roots has these in common with its derivative polynomial. The derivative polynomial can be rationally determined and is of lower degree, hence the original polynomial is not irreducible.

### 3.2 The auxiliary polynomial

**Proposition 3.3.** *For any polynomial with distinct roots one can always find a 'resolvent' meaning a linear rational function V of the roots $r_i$ so that no two values assumed by the resolvent are equal when the roots are permuted*

$$V(r_1..r_n) = a_1 r_1 + a_2 r_2 .. a_n r_n. \tag{31}$$

In the case of an arbitrary polynomial finding $a_i$ does not present any difficulty since the roots are unrestricted, and a choice of $n$ distinct coefficients will suffice. For example, in Eq. (22) we chose the $4^{\text{th}}$ roots of unity, and this can be extended to the general case

$$V(r_1..r_n) = r_1 + \alpha r_2 + ..\alpha^{n-1} r_n, \quad \alpha = e^{2\pi i/n}. \tag{32}$$

A proof of the existence of a resolvent when the coefficients are specified is described in [7]

**Proposition 3.4.** *All roots of a given polynomial can be expressed as rational functions of its resolvent. As a corollary, resolvents differing from each other by a permutation of the roots are all rational functions of one another.*

*Proof.* Consider the irreducible polynomial $p(x)$ with roots $r_1..r_n$

$$p(x) = (x - r_1)(x - r_2)..(x - r_n) = x^n + b_{n-1} x^{n-1} ..+ b_0. \tag{33}$$

Let us construct another polynomial $q(x)$ whose roots are the $(n-1)!$ resolvents formed by permuting all roots of $p(x)$ except $r_1$

$$\begin{aligned} q(x) &= (x - V(r_1, r_2, r_3, r_4..))(x - V(r_1, r_3, r_2, r_4..))(x - V(r_1, r_4, r_3, r_2..))... \\ &= x^n + g_{n-1}(r_1..r_n) x^{n-1} ...+ g_0(r_1..r_n) \end{aligned} \tag{34}$$

Since $q(x)$ is symmetric in $r_2..r_n$, $g_n$ can be rationally expressed in terms of the elementary symmetric polynomials in $n-1$ variables: $\tau_1..\tau_{n-1}$. By Eq. (6) these polynomials are known functions of $r_1$ and the elementary symmetric polynomials in $n$ variables $\tau_i(r_2..r_n) \equiv s_i(r_1, \sigma_1..\sigma_n)$. But $\sigma_1..\sigma_n$ are symmetric in all roots and so are known functions of the coefficients of the original equation, $b_i$. We thus have

$$
\begin{aligned}
g_i(r_1..r_n) &\equiv g'_i(r_1, \tau_1..\tau_{n-1}) \\
&\equiv g'_i(r_1, s_1(r_1, \sigma_1..\sigma_n)..s_{n-1}(r_1, \sigma_1..\sigma_n)) \\
&\equiv g''_i(r_1, b_i) \\
&\equiv g'''_i(r_1)
\end{aligned}
\qquad (35)
$$

We can now re-express Eq. (34) as a polynomial $q(x,r)$ in two variables

$$q(x) \equiv q(x, r_1) = x^n + g'''_{n-1}(r_1)x^{n-1}... + g'''_0(r_1). \qquad (36)$$

By symmetry this can be extended to all roots $r_i$

$$
\begin{aligned}
q(x, r_i) &= x^n + g'''_{n-1}(r_i)x^{n-1}...+ g'''_0(r_i) \\
&= \bigl(x - V(r_i, r_1..r_{j\neq i}..r_n)\bigr)\bigl(x - V(r_i, r_n..r_{j\neq i}..r_1)\bigr).
\end{aligned}
\qquad (37)
$$

where $r_1..r_{j\neq i}..r_n$ denotes a symmetric permutation of the $n-1$ roots $r_1..r_n$ with $r_i$ excluded.

Substituting $V_0 \equiv V(r_1, r_2, r_3..) = a_1 r_1 + a_2 r_2 + ..a_n r_n$ into $q(x,r)$, gives the single-variable polynomial $q(V_0, r)$. Trivially $q(V_0, r_1) = 0$, since $r = r_1$ annuls the first factor in Eq. (37). However $q(V_0, r_2) \neq 0$ since that would mean one of the resolvents with $r_1$ permuted equals resolvent $V_0$ with $r_1$ un-permuted, which contradicts our assumption that all values of the resolvent are different. Therefore $p(r)$ and $q(V_0, r)$ do not share root $r_2$, and by symmetry neither do they share any other root except $r_1$. But this means $r_1$ can be rationally extracted dividing $q(V_0, r)$ by $p(r)$, and is therefore known. We can now repeat the argument with another root $r_2$, and all the other roots, which are therefore also known once one resolvent is given. The roots are therefore all rationally known functions of one resolvent, and therefore by Eq. (31) the resolvents are all rationally known functions of each other, which was to be proved.

**Definition 3.3.** *For a given polynomial with $n$ roots $p(x)$, the auxiliary polynomial $q(x)$ is a monic polynomial whose roots are the $n!$ resolvents of $p(x)$:*

$$q(x) = (x - V(r_1, r_2, ..r_n))(x - V(r_2, r_1, ..r_n))...$$

## 3.3 Reduction of the Galois group by the adjunction of radicals

Finding the roots of a polynomial entails its complete decomposition into linear factors. As demonstrated in Section 3.1, a step in this process is the adjunction of some radical $r = \sqrt[p]{A}$, where $A$ can be any quantity of the type shown in Table 1, which leads to a reduction of the polynomial. Moreover, we need only consider prime roots since extracting any rational root can be decomposed into a sequence of extraction of prime roots and powers. It is well-known that the extraction of a $p$-th root carries an ambiguity associated with multiplication by a $p$-th root of unity. These roots are the solutions of the degree $p-1$ cyclotomic equation presented in Section 4.1, and are therefore solvable in radicals, with the maximum root extracted being $p-1$.

Let us therefore adjoin a root $r_1$ of irreducible polynomial $r^p = A$, to an irreducible factor of the auxiliary polynomial $q(V)$, starting from $p = 2$ and increasing in order, until at some power $p$ the polynomial $q(V)$ reduces. We can write this reduction as the product of a new irreducible polynomial $f(V, r_1)$ and another component

$$q(V) = f(V, r_1) g(V, r_1) \qquad (38)$$

Expanding both sides in powers of $V$ we have

$$q(V) = \sum a_i V^i = f(V, r_1) g(V, r_1) = \sum h_i(r_1) V^i . \qquad (39)$$

This equates known coefficients $a_i$ to the polynomials $h_i(r)$, so that $h_i(r) - a_i = 0$ shares root $r_1$ with $r^p = A$. Therefore $h_i(r) = a_i$ for all $p$ roots of $r^p = A$, and Eq. (38) is satisfied by all $r_i$

$$q(V) = f(V, r_i) g(V, r_i) , \qquad (40)$$

Substituting the $r_i$ into the two variable polynomial $f(V, r)$ generates $p$ polynomials in one variable $f_i(V) \equiv f(V, r_i)$, which from Eq. (40) share all their roots with $q(V)$.

Because of the symmetry between the roots of $r^p = A$, these polynomials are either all identical or all distinct. If they are identical $q(V)$ does not reduce, hence they are distinct, and so the $f_i(V)$ disjointly partition the roots of $q(V)$

$$q(V) = f_0(V) f_1(V) f_2(V) .. f_{p-1}(V) . \qquad (41)$$

Because the RHS of Eq. (41) is symmetric in $r_i$ and therefore a known polynomial, $q(V)$ can admit no other factors. **Thus adjunction of the roots of $r^p = A$ splits $q(V)$ into $p$ irreducible factors,** $f_i(V)$. We now investigate their structure.

Let $V_1$ be a root of $f(V, r_1)$, and say we have found a function $F(V)$ so that $F(V_1)$ is another root. Then the irreducible polynomial $f(V, r_1)$ shares a root with the polynomial $f(F(V), r_1)$. It therefore shares all roots, and thus divides it

$$f(F(V), r_1) = f(V, r_1) g(V, r_1), \tag{42}$$

where $g(V, r)$ is the other factor. Repeating the arguments which lead from Eq. (38) to Eq. (40) this must hold for all roots $r_i$

$$f(F(V), r_i) = f(V, r_i) g(V, r_i). \tag{43}$$

Therefore if $V_i$ is a root of $f(V, r_i)$ then $F(V_i)$ is another root. Putting this another way, if the function $F(V)$ relates two roots of $f_i(V)$, it also relates two roots of $f_j(V)$.

Let the $V_i$ and $V_j$ be related by the permutation $S$ in Eq. (31). Then the values of $F(V_i)$ and $F(V_j)$ are related by the same permutation, as are therefore the roots of $f_i(V)$ and $f_j(V)$. This is summarized in the table below, where $S = (1,2)$ for clarity

$$
\begin{array}{ll}
f_i(V) & f_j(V) \\
\\
V_i \equiv V(r_1, r_2, ..r_n) \xrightarrow{S} V_j \equiv V(r_2, r_1 ..r_n) \\
F_2(V(r_1, r_2, ..r_n)) \xrightarrow{S} F_2(V(r_2, r_1, ..r_n)) \\
F_3(V(r_1, r_2, ..r_n)) \xrightarrow{S} F_3(V(r_2, r_1, ..r_n)) \\
......... \qquad\qquad\qquad ......... \\
F_k(V(r_1, r_2, ..r_n)) \xrightarrow{S} F_k(V(r_2, r_1, ..r_n))
\end{array}
\tag{44}
$$

Each column lists the $k$ roots of $f(V)$, $V_k = F_k(V)$. Since $f_i(V)$ is irreducible, functions invariant with respect to all $V_k$ are known functions of the coefficients of $f_i(V)$ after root $r_i$ is adjoined. On the other hand, from Eq. (31), the set $V_k$ is in a one-to-one correspondence to permutations of the roots of $p(x)$. It follows (see Section 4) that these permutations correspond to the Galois group of $f_i(V)$, which is thus a subgroup of the Galois group of $q(V)$. Because of the symmetry between the roots of $r^p = A$ this subgroup is identical for all $i$, and therefore the columns in Eq. (44) correspond to its presentations. However Eq. (44) shows that the presentation are *also* related by a single permutation $S$ (as in Eq. (16)). Since both left and right cosets are therefore identical it follows that the Galois subgroup is normal, and the prime number of presentations means the quotient group is cyclic of prime order. **The adjunction of the roots of $r^p = A$ to $q(V)$ reduces its Galois group to a normal subgroup, with a cyclic quotient group of prime order.** We thus arrive at the solution of Section 2.7, which is therefore unique.

## 4. Computing the Galois group

In Section 3.2 we found a function $q(V, x)$ such that if $q(V_0, r_1) = 0$, $r_1$ can be obtained by polynomial division, $q(V_0, x) / p(x)$. Let us write the result of this division as $r_1 \equiv F(V_0)$, where $F(V)$ is known. Replacing resolvent $V_0$ by resolvent $V_1 \equiv a_1 r_2 + a_2 r_1 + ..a_n r_n$, with $r_2$ and $r_1$ swapped, gives $q(V_1, r_2) = 0$. The inverse of this relationship is the same as above, $r_2 = F(V_1)$. Thus if a function $F(V)$ returns a root of $p(x)$ for some resolvent, it will return a root of $p(x)$ for all resolvents.

Given this, assume we have found a set of functions $F_i(V)$ each returning a root $r_i$ when acting on the resolvent $V_0$, $r_i \equiv F_i(V_0)$. Assume also we have found a set of resolvents $V_0..V_m$ all belonging to the same irreducible factor of the auxiliary equation. Consider the following table

$$\begin{array}{cccc} F_1(V_0) & F_2(V_0) & F_3(V_0) & ... & F_n(V_0) \\ F_1(V_1) & F_2(V_1) & F_3(V_1) & ... & F_n(V_1) \\ ... & ... & ... & & ... \\ F_1(V_m) & F_2(V_m) & F_3(V_m) & ... & F_n(V_m) \end{array} \tag{45}$$

Since all functions in the table return a root of $p(x)$, with the roots arranged as $r_1..r_n$ in the top row, the table actually represents a set of permutations of the roots $r_1..r_n$. We now claim:

**Proposition 4.1** *The set of permutations in Eq. (45) forms the Galois group of $p(x)$, with all functions invariant under these permutations being known, and vice versa, all known functions being invariant under this set of permutations.*

To prove this, consider a function of the roots $f(r_1, r_2..r_n)$ invariant under the permutations of Eq. (45). We can write this invariance as

$$f(F_1(V_0), F_2(V_0)..F_n(V_0)) = f(F_1(V_1), F_2(V_1)..F_n(V_1)) = ... \tag{46}$$

thus $f(r_1, r_2..r_n) = f(F_1(V), F_2(V),..F_n(V))$ is invariant for all $V_i$. It is therefore a symmetric function in the roots $V_i$, and is therefore known.

Conversely, let us write a known function of the roots $f(r_1, r_2..r_n)$ as

$f(r_1, r_2..r_n) = f(F_1(V_0), F_2(V_0)..F_n(V_0)) \equiv g(V_0)$, where $g(V)$ is therefore also known. If $f(r_1, r_2..r_n)$ is not invariant under all permutations of the above table, $g(V)$ is not invariant under all $V_i$. Since this means $g(V) - \text{const} = 0$ for only some $V_i$, we have the contradiction of a known polynomial which does not share all its roots with an irreducible polynomial.

Before turning to applications we note the following. Firstly, one generally does not have to construct the functions $F_i(V)$. Below we calculate the Galois group of a quintic by finding an irreducible factor of its auxiliary polynomial, and establishing the permutations directly. Secondly, the resolvents $V_i$ in Eq. (45) must belong to the same irreducible factor of the auxiliary equation. If on the other hand all $n!$ resolvents of the auxiliary equation are used, and for some choice of functions $F_i(V)$ the group of permutations is reduced to a subgroup of $S_n$, then clearly the Galois group is contained in this group. We use this to calculate the Galois group of the cyclotomic equation.

## 4.1 Galois group of the cyclotomic equation of prime order

The roots of all cyclotomic equations of degree $p$

$$(\alpha^p - 1)/(\alpha - 1) \equiv \alpha^{p-1} + \alpha^{p-2} + ..1 = 0, \qquad (47)$$

can be expressed in radicals in terms of the roots of cyclotomic equations of prime degree, hence only these will be considered.

It is well known that for every prime $p$ there exists a primitive integer $g$ such that all integers between 1 and $p-1$ can be expressed as powers of $g$ modulus $p$, $g^n \mod p = i$. This means that if $\alpha$ is a $p$-th root of unity, the sequence $\alpha, \alpha^g, (\alpha^g)^g, ((\alpha^g)^g)^g$, generates all roots of unity.

Example: For $p = 11$, a primitive index is $g = 2$, and so $\alpha_i = \alpha^{g^i}$ gives

$$\alpha_i = \{\alpha, \alpha^2, \alpha^4, \alpha^8, \alpha^5, \alpha^{10}, \alpha^9, \alpha^7, \alpha^3, \alpha^6\}. \tag{48}$$

Let us define $f(r) = r^g$ where $g$ is the primitive index, choose a root $\alpha_1$ and a resolvent $V_0$, and find a function $F(V)$ such that $\alpha_1 = F(V_0)$. Then the entire table of Eq. (45) can be generated from just $F(V)$ and $f(r)$

$$\begin{aligned} F(V_0) &\equiv \alpha_1 & f(F(V_0)) &= \alpha_2 & f^2(F(V_0)) &= \alpha_3 & \ldots & f^{n-1}(F(V_0)) &= \alpha_{n-1} \\ F(V_1) &\equiv \alpha_i & f(F(V_1)) &= \alpha_{i+1} & f^2(F(V_1)) &= \alpha_{i+2} & \ldots & f^{n-1}(F(V_1)) &= \alpha_{i+n-1} \\ &\ldots & &\ldots & &\ldots & & &\ldots \\ F(V_m) &\equiv \alpha_k & f(F(V_m)) &= \alpha_{k+1} & f^2(F(V_m)) &= \alpha_{k+2} & \ldots & f^{n-1}(F(V_m)) &= \alpha_{k+n-1} \end{aligned} \tag{49}$$

where $f(f(r)) \equiv f^2(r)$ etc., and $\alpha_1, \alpha_i..\alpha_k$ in the first column denote whatever roots our chosen function $F(V)$ generates for the given argument. We see that the Galois group of the cyclotomic equation of prime degree is the group of all cyclic permutations $Z_{p-1}$. This outcome is independent of the choice of functions $F_i(V) = f^i(F(V))$ in the table, since a different choice simply changes the starting permutation in the top row.

Example: Let $\alpha$ and $\beta$ be the $17^{th}$ and the (known) $16^{th}$-roots of unity respectively. Let

$$t \equiv \alpha + \beta\alpha^2 + \beta^2\alpha^3 .. + \beta^{15}\alpha^{16}. \tag{50}$$

Cyclically permuting the roots $\alpha$ is equivalent to multiplying $t$ by some power of $\beta$, a 16-th root of unity, and so $t^{16}$ is invariant under cyclic permutations. Therefore $t^{16}$ is symmetric under the Galois group of the equation, and is known; indeed expanding $t^{16}$ we find it is a function of $\beta$ alone. Since there is nothing special about the particular permutation of $\alpha^i$ in Eq. (50) *all* such combinations are known, from which $\alpha$ can be extracted by a solution of 16 simultaneous equations. Since determining $t$ involves extracting the 16-th root of unity, which can be expressed in square roots, the regular 17-gon can be drawn with ruler and compass.

## 4.2 Galois group of a quintic

The general procedure is as follows. The Galois group of a solvable quintic is a subgroup of $F_{20}$ which, as we saw, consists of 4 presentations of the cyclic group $Z_5$ related by the cycle (2,3,5,4). Choosing a suitable resolvent $t$, we compute the polynomial

$$P(x) = (x-t_1)(x-t_2)..(x-t_{20}), \tag{51}$$

whose roots are the images of $t$ under the 20 substitutions of $F_{20}$. If this polynomial has integer coefficients the quintic is solvable. For this it is sufficient (by Lagrange's construction) to check that $P(x)$ is (complex) integer valued for 20 different (complex) integers $x_i$. The Galois group can be further narrowed down to any of the solvable subgroups of $F_{20}$ by restricting the polynomial $P(x)$ to just the permutations of this subgroup.

Example: Consider DeMoivre's quintic $x^5 + 5ax^3 + 5a^2 x + b = 0$. Choosing $a = -2, b = 1$ we evaluate the roots numerically and construct the following resolvent:

$$t = (2.705..) + i(0.050..) - (2.674..) - i(-1.702..) + 0(1.621..) \tag{52}$$

where $t$ assumes different values under all permutation of the roots $r_i$. Choosing the 25 complex integers $x_i = c_i + id_i$ with $-2 \leq c_i, d_i \leq 2$, we calculate the product Eq. (51). This can be done using standard 15-digit computation in Excel with almost no rounding errors. A product of 20 resolvents $\prod(x_i - t_p)$ gives a clean integer on the twentieth multiplication, and $P(x_i)$ is integer valued for all choices of $x_i$. Therefore $P(x_i)$ has integer coefficients, and $x^5 + 5x^3 + 5x + 1 = 0$ is solvable. Similarly, other choices of $a, b$ lead to solvable polynomials. Therefore a solution in radicals can be found by following the general procedure of Section 2.7, or, in this case, more simply using the substitution $x = y - a/y$, to give the five roots $r_m$

$$r_m = \alpha^m u_1^{1/5} + \alpha^{5-m} u_2^{1/5}, \quad u_1, u_2 = \left(-b \pm \sqrt{b^2 + 4a^5}\right)/2, \tag{53}$$

where $\alpha$ is a 5$^{th}$ root of unity. If $b^2 + 4a^5 > 0$ there is one real solution and two complex conjugate pairs, in the contrary case $u_1^* = u_2$, and all five solutions are real. In this case it is more convenient to re-express Eq. (53) in the form

$$r_m = 2(-a)^{1/2} \cos\left(\frac{1}{5}\left[2\pi n + \cos^{-1}(-b/2(-a)^{5/2})\right]\right). \tag{54}$$

Example: Consider the irreducible binomial $x^5 = A$. From Eq. (41) its Galois group is $Z_5$ once the 5$^{th}$ roots of unity are adjoined. The Galois group of the 5$^{th}$ roots of unity, from Eq. (49) is $Z_4$, so the Galois group of the binomial is composed of 4 presentations of $Z_5$ linked by permutations of $Z_4$, this corresponds to the Frobenius group $F_{20}$.

## *5. Conclusion*

A method for solving polynomial equations in radicals based on simple group theory methods has been presented and applied to solve general polynomials of degree 4 and below, and specific quintics. Polynomials with degree $\geq 5$ with arbitrary coefficients are shown to be insoluble by this method. Galois theory was used to demonstrate that no other solution in radicals is possible. A method for calculating the Galois group was presented and used to solve particular polynomial equations.

1. L. Novotny, Strong coupling, energy splitting, and level crossings: A classical perspective, *Am. J. Phys.*, *78*(2010), 1199-1202.

2. G. Toth, *Groups Minimal Immersion of Spheres, and Moduli*, Springer, New York, 2002, p. 66.

3. J. Stillwell, *Galois theory for beginners*, *Am. Math. Monthly*. *101*(1994) 22-27.

4. E. Galois, Oeuvres mathematiques d'Evariste Galois, *J. Math. Pur. Appl. 11*(1844) 381-444.

5. E. Artin, *Galois Theory*, Dover, New York, 1998.

6. I. Stewart, *Galois* Theory, Chapman & Hall, London, 1976.

7. J. Bewersdorff, *Galois Theory for Beginners*, AMS, Providence, 2006.

8. H. M. Edwards, *Galois Theory*, Springer-Verlag, New York, 1984.

9. L. Sylow and S. Lie, *Oeuvres Completes De Niels Henrik Abel*, Kessinger Publishing, Whitefish MT, 2010.

10. D.S. Dummit, *Math. Comp*. *195*(1991) 387-401.

11. B. Spearman and K. S. Williams, Characterization of solvable quintics $x^5 + ax + b$, *Am. Math. Monthly*., *101*(1994) 986–992.

**Figure Captions**

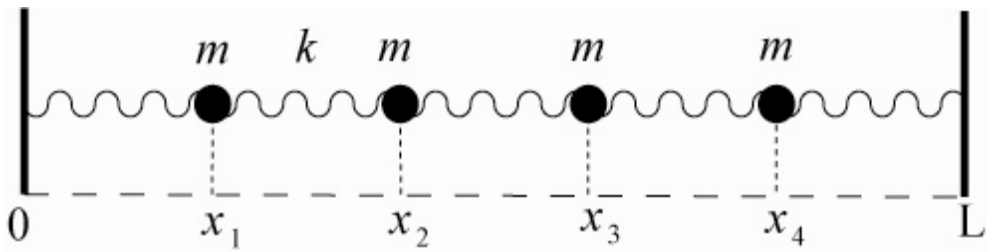

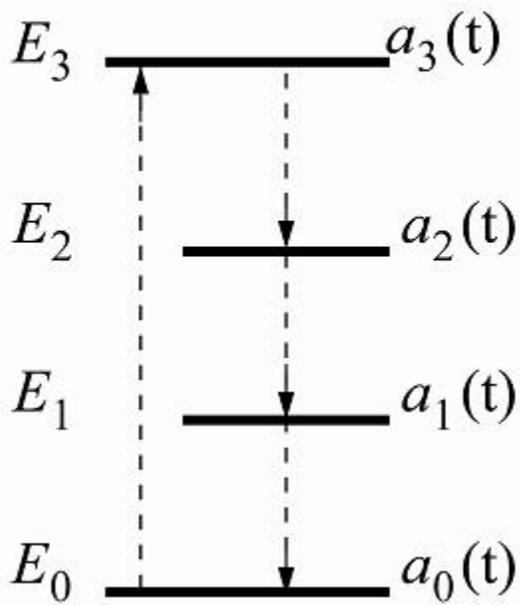

Figure 1. Oscillatory systems with frequencies and normal modes determined by the roots of polynomial equations. (a) A set of masses connected by springs and (b), a quantum system approximated by a set of spaced quantum levels.

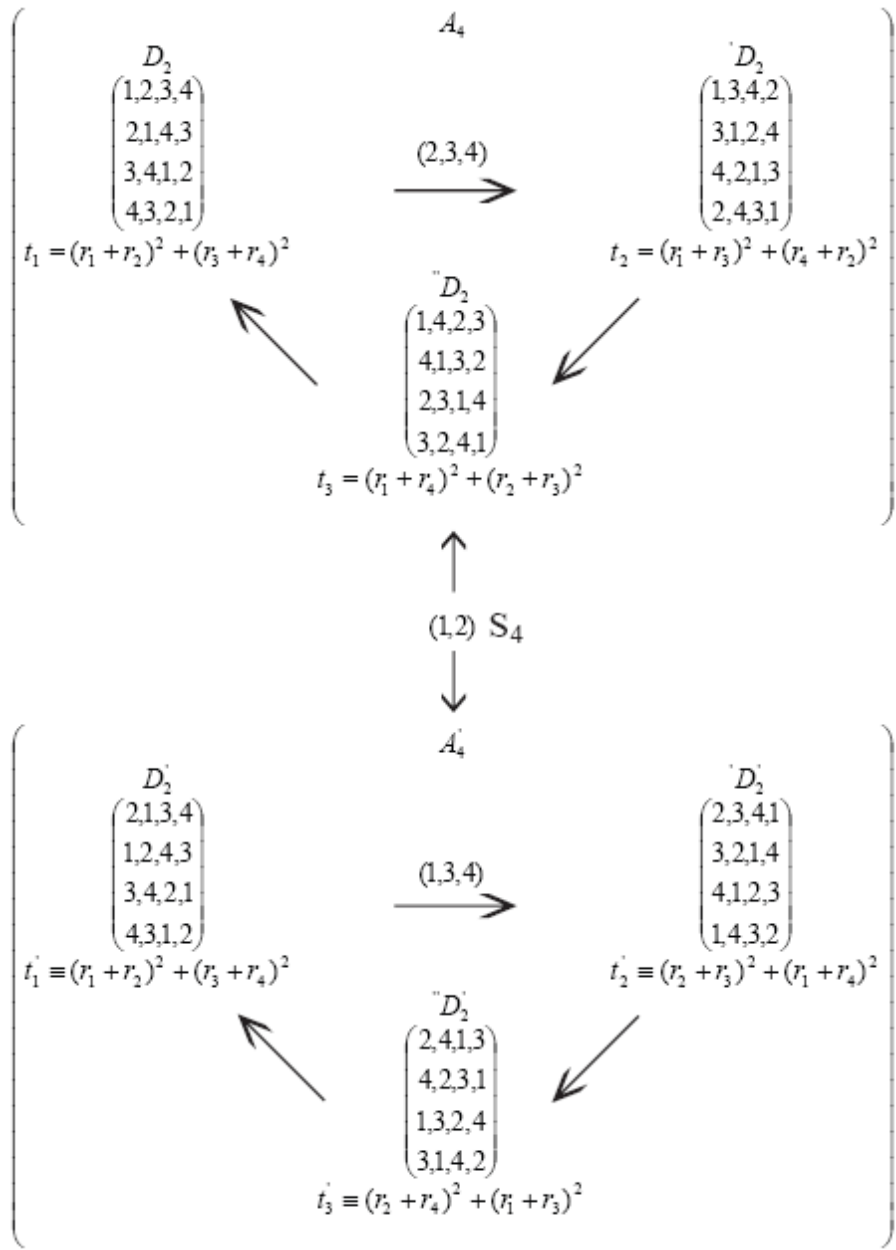

Figure 2. Structure of the Galois group of the general quartic $S_4$, showing that it is solvable in radicals. The 24 elements of the group can be composed into a series of subgroups $D_2 \triangleleft A_4 \triangleleft S_4$, with presentations linked by a single permutation, the generator of the quotient group. The diagram also shows the polynomials $t_i$ invariant under each presentation of $D_2$.

35